\documentclass[letterpaper, 10 pt, conference]{ieeeconf}  
\usepackage{cite}
\usepackage{amsmath,amssymb,amsfonts}
\usepackage{algorithmic}
\usepackage{graphicx}
\usepackage{subfigure} 
\usepackage{textcomp}
\usepackage{xcolor}
\usepackage{algorithm}
\usepackage{algorithmic}

\IEEEoverridecommandlockouts                              

\title{\LARGE \bf
Differentially Private Distributed Mismatch Tracking Algorithm for Constraint-Coupled Resource Allocation Problems
}

\author{Wenwen Wu$^{1}$, Shanying Zhu$^{1}$, Shuai Liu$^{2}$ and Xinping Guan$^{1}$
\thanks{*This work was supported in part by the NSF of China under the Grants 61922058, 62173225 and the Key Program of the National Natural Science Foundation of China under Grant 62133008.}
\thanks{$^{1}$Wenwen Wu, Shanying Zhu and Xinping Guan are with the Department of Automation, Shanghai Jiao Tong University, Shanghai 200240, China; Key Laboratory of System Control and Information Processing, Ministry of Education of China,Shanghai 200240, China, and also Shanghai Engineering Research Center of Intelligent Control and Management, Shanghai 200240, China.}%
\thanks{$^{2}$Shuai Liu is with the School of Control Science and Engineering, Shandong University, Jinan,China.}%
\thanks {Final version of the article available at https://doi.org/10.1109/CDC 51059. 2022.9993173  (please cite as \cite{DMAC}).}%
}

\begin{document}

\maketitle
\thispagestyle{empty}
\pagestyle{empty}

\begin{abstract}

This paper considers privacy-concerned distributed constraint-coupled resource allocation problems over an undirected network, where each agent holds a private cost function and obtains the solution via only local communication. With privacy concerns, we mask the exchanged information with independent Laplace noise against attackers with potential access to all network communications. We propose a differentially private distributed mismatch tracking algorithm (diff-DMAC) to achieve cost-optimal distribution of resources while preserving privacy. Adopting constant stepsizes, the linear convergence property of diff-DMAC in mean square is established under the standard assumptions of Lipschitz smoothness and strong convexity. Moreover, it is theoretically proven that the proposed algorithm is $\epsilon$-differentially private. And we also show the trade-off between convergence accuracy and privacy level. Finally, a numerical example is provided for verification. 

\end{abstract}

\section{Introduction}

Distributed resource allocation problem (DRAP) has received much attention due to its wide applications in various domains, including smart grids \cite{economic_dispatch}, wireless networks\cite{wireless_networks} and robot networks \cite{robot_task_allocation}. The goal of DRAP is to achieve the  cost-optimal distribution of limited resources among users to meet their demands, local constraints, and possibly certain coupled global constraints, see \cite{Tracking-ADMM} for example. Conventional centralized approaches are subject to performance limitations of the central entity, such as a single point of failure, limited scalability. Moreover, it may raise privacy concerns when the central agent is not reliable enough. Alternatively, distributed approaches which operate with only local information have better robustness and scalability, especially for large-scale systems \cite{DuSPA,DDGT,distributed_ADMM_like}. 

To implement algorithms in a distributed manner, the information exchange between agents in the network is unavoidable, which can raise concerns about privacy disclosure. For example, in the supply–demand optimization of the power grid, attackers can infer users' life patterns from certain published information that is computed based on the demand from users\cite{2017_TAC_diff_constrained}. While the aforementioned works consider an ambitious suite of topics under various constraints imposed by real-world applications, privacy is an aspect generally absent in their treatments.

Recently, privacy preservation becomes an increasingly critical issue for distributed systems. Encryption-based algorithms are proposed for distributed optimization problems \cite{encryption1}, \cite{encryption2} and proved to be effective. However, the encryption technique leads to a high computational complexity which seems undesirable in large-scale networks. To preserve the privacy of the agents’ costs, an asynchronous heterogeneous-stepsize algorithm is proposed in \cite{asynchronous_heterogeneous_stepsize_algorithm}, but the method assumes
that the information of communication topology is unavailable to attackers. 

Another thread to address distributed optimization problems is to adopt  differential privacy technique which is robust to arbitrary auxiliary information exposed to attackers, including communication topology information, thus well suited for multi-agents scenarios. Motivated by the privacy concerns in EV charging networks, the work in \cite{2017_TAC_diff_constrained} presents a differentially private distributed algorithm to solve distributed optimization problem, but an extra central agent is needed. In \cite{2018_TSIPN_decay_stepsize}, the authors design a completely distributed algorithm that guarantees differential privacy by perturbing the states with Laplace noise. It requires the stepsize to be decaying to guarantee the convergence, resulting in a low convergence rate. Adopting the gradient-tracking scheme, a differentially private distributed algorithm is proposed in \cite{TAC_Ding} which has a linear convergence rate with constant stepsizes. However, the proposed algorithm in \cite{TAC_Ding} is only applicable to unconstrained problems. For solving the economic dispatch problem where both global and local constraints are considered, a privacy-preserving distributed optimization algorithm is proposed in \cite{2020_TII} for quadratic cost functions and its effectiveness in preserving privacy is guaranteed through qualitative analysis. The extension to general DRAPs is performed in \cite{2021_TSIPN_diff_DRADT} for strongly convex and Lipschitz smooth cost functions. Additionally, quantitative analysis of privacy is provided. However, the work only considers scalar cost functions and the coupling within constraints is ignored.

In this paper, a differentially private distributed mismatch tracking algorithm (diff-DMAC) is proposed to solve constraint-coupled DRAPs with individual constraints and a global linear coupling constraint while preserving the privacy of the local cost functions. To meet the global coupled-constraint, a mismatch tracking step is introduced to obtain the supply-demand mismatch in the network. And the tracked global mismatch is then implemented for error compensation. With a constant stepsize, diff-DMAC has a provable linear convergence rate for strongly convex and Lipschitz smooth cost functions. We further present a rigorous analysis of the algorithm's differential privacy and thus provide a stronger privacy guarantee. Finally, its convergence properties and effectiveness in privacy-preserving are numerically validated.

The rest of the paper is organized as follows. Section \ref{Problem Formulation} formulates the constraint-coupled DRAP. The proposed algorithm is developed in Section \ref{Algorithm Development}. And the theoretical analyses about its convergence and privacy are given in Section \ref{convergence analysis} and \ref{Differential Privacy}, respectively. Then the algorithm is numerically tested in Section \ref{Numerical Experiments}. Finally, Section \ref{Conclusion} concludes the paper. 

\emph{Notation:} Vectors default to columns if not otherwise specified. Bold letter $\mathbf{x}\in \mathbb{R}^{n \times p}$ is defined as $\mathbf{x}=[x_1^\intercal, \cdots, x_n^\intercal]^\intercal$. The Kronecker product is denoted by $\otimes$. Let $\mathbf{1}_n$ be the $n$-dimension vector with all one entries. For vectors, we use $\Vert \cdot \Vert$ to denote the 2-norm. And for matrices, $\Vert \cdot \Vert$ denotes the spectral norm. $\underline{\lambda}(X)$ denotes the minimum eigenvalue of matrix $X$. We use $blkdiag(X_1,\cdots,X_n)$ to refer to the block-diagonal matrix
with $X_1,\cdots,X_n$ as blocks. For a random variable $x \in \mathbb{R}$, $Lap(\theta)$ denotes the Laplace distribution with pdf $f=(x \vert \theta)=\frac{1}{2\theta}e^{-\frac{\vert x \vert}{\theta}}$, where $\theta>0$. If $x \sim Lap(\theta)$, we have $\mathbb{E} [\vert x \vert ]=\theta$, $\mathbb{E} [x^2]=2\theta^2$. $\mathbb{P}(\cdot)$ is used to denote the probability. $\mathbb{B}(S)$ denotes the set of Borel subsets of topological space $S$.

\section{Problem Formulation}
\label{Problem Formulation}
We consider the following general constraint-coupled distributed resource allocation problem
\begin{align}
	\min _{\boldsymbol{x}\in \mathbb {R}^{np}}f(\boldsymbol{x})&=\sum _{i=1}^{n}f_i(x_i) \notag \\
	s.t. \sum _{i=1}^{n} A_i x_i&=\sum _{i=1}^{n}d_i,  \ \,x_i\in \mathcal {X}_i,\forall i,
	\label{DRAP}
\end{align}
where  $f_i:\mathbb{R}^p \rightarrow \mathbb{R}$ is the agent $i$'s private cost function,  $x_i \in \mathbb{R}^p$ is  the  decision vector of agent $i$ and $d_i$ denotes the local resource demand. $\mathcal{X}_i$ is a local convex and closed set which encodes local constraints of agent $i$. $A_i \in \mathbb{R}^{m\times p} \ (p\geq m)$ is the coupling matrix with full row rank and $A_i A_i^\intercal$ is invertible, i.e., $\underline{\lambda}(A_i A_i^\intercal)>0$.

\begin{figure}[thpb]
	\centering
	\includegraphics[width=3in]{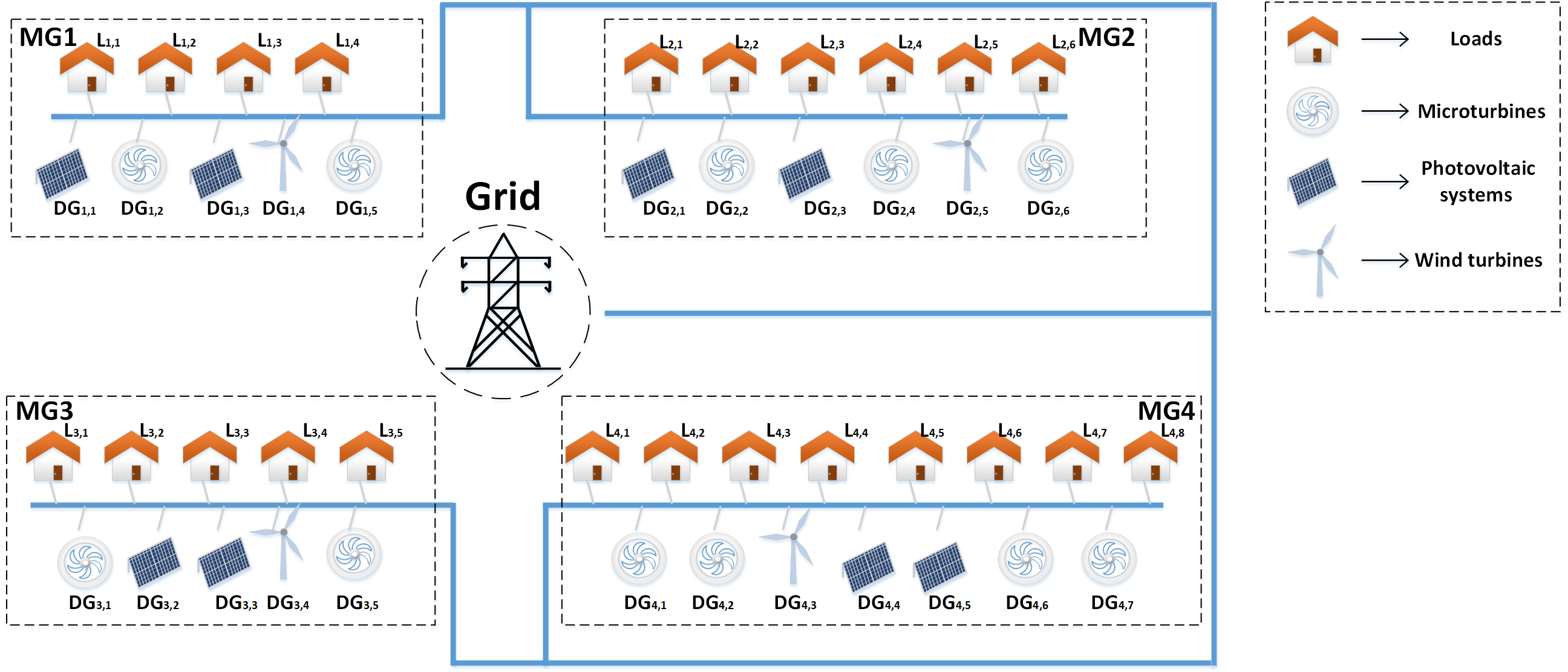}
	\caption{Multi-microgrid systems with distributed generators.}
	\label{fig}
\end{figure}

Many practical problems take the form of problem (\ref{DRAP}), e.g.,\cite{economic_dispatch,multi_microgrid_systems}. One typical example is the economic dispatch of multi-microgrid (multi-MG) systems as shown in Fig. \ref{fig}, where both microturbines and non-dispatchable distributed generators are involved to provide energy to meet the loads, and the outputs of non-dispatchable distributed generators can fluctuate significantly \cite{multi_microgrid_systems}. When the supply-demand balance in the MG cannot be maintained, the
coordination among MGs is needed and both within-MG, between-MG optimization problems should be considered to improve systems' reliability \cite{MG_connected}.

\emph{Assumption 1:} Each cost function $f_i:\mathbb {R}^{p} \rightarrow \mathbb {R}$ is strongly convex and Lipschitz smooth, i.e.,
\begin{equation}
	\Vert \nabla f_i(x)-\nabla f_i(x')\Vert \leq L_i\Vert x-x'\Vert
\end{equation}
and
\begin{equation}
	(x-x')^\intercal \left(\nabla f_i(x)-\nabla f_i(x')\right)\geq \varphi _i \Vert x-x'\Vert^2
\end{equation}
$\forall x,x'\in \mathbb {R}^p$, where $\varphi_i >0$ and $L_i < \infty$ are the strong convexity and Lipschitz constants, respectively. Define $\underline{\varphi}=\min_i\{\varphi_i\}$ and $\overline{L}=\max_i\{L_i\}$.

The communication network over which agents exchange information can be represented by an undirected graph $\mathcal{G}=\left( \mathcal{N},\mathcal{E} \right) $ where $\mathcal{N}=\{1,\cdots,n\}$ is the set of agents and $\mathcal{E}\subseteq \mathcal{N}\times \mathcal{N}$ denotes the set of edges, accompanied with a nonnegative weighted matrix $\mathcal{W}$. For any $ i,j \in \mathcal{N}$ in the network, $w_{ij}>0$ denotes agent $j$ can exchange information with agent $i$. The collection of all individual agents that agent $i$ can communicate with is defined as its neighbors set $\mathcal{N}_i$.

\emph{Assumption 2:} The graph $\mathcal {G}$ is undirected and connected and the weight matrix $\mathcal{W} \in \mathbb{R}^{n \times n}$ is doubly stochastic, i.e., $\mathcal{W}\mathbf{1}_n=\mathbf{1}_n$ and $\mathbf{1}_n^{\intercal}\mathcal{W}=\mathbf{1}_n^{\intercal}$.

Under Assumption 2, from the Perron–Frobenius theorem, we have that 
	$\bar{\lambda}=\Vert (\mathcal{W}-\frac{\mathbf{1}^{\intercal}\mathbf{1}}{n}) \otimes I_m \Vert<1$.

Assumptions 1-2 are standard when solving related problems.

In this problem, agents want to achieve the optimum while preserving the privacy of local cost functions which can be commercially sensitive in practice. To preserve privacy of functions, we add noise $\boldsymbol{\zeta }$ to mask  the exchanged information $\mathbf{x}(k)$ at each time instant as $\mathbf{z}(k)=\mathbf{x}(k)+\boldsymbol{\zeta }(k)$. Let $\mathcal{F}$ be the function set $\{f_i\}_{i=1}^n$, the update of $\mathbf{x}$ can be expressed in impact form 
\begin{align*}
	\mathbf{x}(k+1)=g(\mathbf{x}(k),\boldsymbol{\zeta }(k),\mathcal{W},\mathcal{F}),
\end{align*}
where $g(\cdot)$ denotes the update rule. Therefore, exchanged information $\mathbf{z}=\{\mathbf{z}(k)\}_{k=0}^\infty $ can be determined if the initialization $\mathbf{x}(0)$, the added noises $\boldsymbol{\zeta }=\{\boldsymbol{\zeta }(k) \}_{k=0}^\infty$, the communication topology $\mathcal{W}$ and the function set $\mathcal{F}$ are known. To provide theoretical guarantees on the privacy for their cost functions, the concept of differential privacy is introduced, see \cite{2021_TSIPN_diff_DRADT}.

\emph{Definition 1:} Given $\delta>0$ and two function sets $\mathcal{F}^{(1)}=\{f_i^{(1)}\}_{i=1}^n$ and $\mathcal{F}^{(2)}=\{f_i^{(2)}\}_{i=1}^n$, the two subsets are $\delta$-adjacent if there exists $\delta' \in (-\delta,\delta)$ and some $i_0$ such that 
\begin{equation*}
	f_i^{(1)}=f_i^{(2)}, \forall i\ne i_0; \nabla f_{i_0}^{(1)}(x)=\nabla f_{i_0}^{(2)}(x+\delta ^{\prime }), \forall x\in \mathcal {X}_{i_0}^{(1)}, 
\end{equation*}
where $\mathcal{X}_{i_0}^{(l)}$ is the domain of $f^{(l)}_{i_0},l=1,2$ and they satisfy that $\mathcal {X}_{i_0}^{(2)}=\lbrace x+\delta ^{\prime }|x\in \mathcal {X}_{i_0}^{(1)}\rbrace$.

\emph{Definition 2:} Given $\delta, \epsilon >0$, for any two $\delta$-adjacent function sets $\mathcal{F}^{(1)},\mathcal{F}^{(2)}$ and any observation $\mathcal {O}\subset \mathbb {B}((\mathbb {R}^{nm})^\mathbb {N})$, the process keeps $\epsilon$-differentially private if
\begin{align*}
	\mathbb {P}\lbrace \boldsymbol{\zeta }\in \Omega |Z_{T^{(1)}}(\boldsymbol{\zeta })\in \mathcal {O}\rbrace
	\leq e^{\epsilon }\mathbb {P}\lbrace \boldsymbol{\zeta }\in \Omega |Z_{T^{(2)}}(\boldsymbol{\zeta })\in \mathcal {O}\rbrace,
\end{align*}
where we define $T^{(l)}=\{x(0), \mathcal{W}, \mathcal {F}^{(l)}\}$, $Z_{T^{(1)}}(\boldsymbol{\zeta })=\mathbf{z}$ and the observation $\mathcal {O}$ encodes all the information collected by the potential attacker. The constant $\epsilon$ indicates the level of privacy: smaller $\epsilon$ implies higher level of privacy.

The attacker is assumed to know all auxiliary information, including the exchanged information between agents, communication topology, etc., except the parameter $\delta'$. Intuitively, differential privacy guarantees that the two function sets can not be distinguished via the released information and thus avoids the leakage of the optimal operating point. 
\section{Algorithm Development}
\label{Algorithm Development}
In this section, a distributed algorithm is designed to solve constraint-coupled DRAPs with privacy concerns. The following Proposition is presented to give the optimal conditions for constraint-coupled DRAP.

\emph{Proposition 1:} $x^\star=[{x_1^\star}^\intercal, \cdots,{x_n^\star}^\intercal] \in \mathbb{R}^{n \times p}$ is the optimal solution of (\ref{DRAP}) if it satisfies \\
\romannumeral1) $\ \sum _{i=1}^n A_i x_i^\star = \sum _{i=1}^n d_i$, \\
\romannumeral2) there exists $\mu^\star \in \mathbb{R}^m$ such that
\begin{align}
	x^\star= \arg \min _{\mathbf{z} \in \mathcal{X}} \{f(\mathbf{z})- (\mathbf{1}_n \otimes \mu^\star)^\intercal \mathbf{A}\mathbf{z} \}, \label{Propositon1_x}
\end{align}
where  $\mathcal{X}$ is defined as the Cartesian product $\mathcal{X} = \mathcal{X}_1 \times \cdots \times \mathcal{X}_n$ and the matrix  $\mathbf{A}= blkdiag(A_1,\cdots,A_n)$. 

\emph{Proof:} \romannumeral1)  is one of the feasible conditions of (\ref{DRAP}) which must be satisfied. Problem (\ref{Propositon1_x}) can be decomposed into $n$ subproblems $x_i^\star=\arg \min _{z\in \mathcal {X}_i}\lbrace f_i(z)-{\mu ^\star}^{\intercal} A_i z \rbrace, \forall i=1, \cdots,n$. It follows that 
\begin{align*}
	f_i(x_i^\star)-{\mu ^\star}^{\intercal} A_i x_i^\star \leq f_i(x_i)-{\mu ^\star}^{\intercal} A_i x_i, \ \forall x_i \in \mathcal{X}_i.
\end{align*}
Summing the above inequality from $i=1$ to $n$ yields
\begin{align*}
	\sum \limits_{i=1}^n 	f_i(x_i^\star) \leq \sum \limits_{i=1}^n  f_i(x_i)-{\mu ^\star}^{\intercal} ( \sum \limits_{i=1}^n A_ix_i - \sum \limits_{i=1}^n A_ix_i^\star).
\end{align*}
If $\mathbf{x}^\star$ satisfies condition \romannumeral1), $\sum _{i=1}^n A_i x_i - \sum_{i=1}^n A_i x_i^\star =0$ holds for any feasible solution $\mathbf{x}$ of (\ref{DRAP}). Thus, we have $f(\mathbf{x}^\star) \leq f(\mathbf{x})$ which completes the proof. $\hfill\blacksquare$

In the following, we will use Proposition 1 to develop a distributed algorithm for solving the DRAP (\ref{DRAP}). Based on the condition \romannumeral2) of Proposition 1, the $\mathbf{x}$-update is designed and the consensus protocol \cite{consensus_protocal} is adopted to drive $\mu_i$ in each agent to a same value
\begin{align}
	&\boldsymbol{\mu}(k+1)=\mathbf{W}\boldsymbol{\mu}(k), \label{compact_mu_1} \\
	&\mathbf{x}(k+1)=\arg \min _{\mathbf{z}\in \mathcal {X}}\lbrace  f(\mathbf{z})-\boldsymbol{\mu}^{\intercal}(k+1) \mathbf{A}\mathbf{z} \rbrace, \label{compact_x_1}
\end{align}
where we define $\mathbf{W}=\mathcal{W}\otimes I_m$.

Inspired by the gradient-tracking scheme \cite{gradient_tracking} , a tracking step is introduced here to track the global supply-demand mismatch as follows 
\begin{align}
\mathbf{y}(k+1)=\mathbf{W}\mathbf{y}(k)+\mathbf{A}\mathbf{x}(k+1)-\mathbf{A}\mathbf{x}(k). \label{deviation_tracking}
\end{align}
Under Assumption 2, if the initialization condition $(\mathbf{1}_n \otimes I_m)^\intercal (\mathbf{A}\mathbf{x}(0)-\mathbf{d})=(\mathbf{1}_n \otimes I_m)^\intercal \mathbf{y}(0)$ holds, iteration (\ref{deviation_tracking}) can infer that
\begin{align}
	\sum _{i=1}^n y_i(k)=\sum _{i=1}^n A_ix_i(k) - \sum _{i=1}^n d_i, \ \forall k>0.  \label{tracking_expression}
\end{align}
To guarantee the condition \romannumeral1) of Proposition 1 holds, the tracked mismatch is used as a compensation term for $\boldsymbol{\mu}$-update (\ref{compact_mu_1}). Thus, the proposed algorithm is obtained
\begin{align}
	&\boldsymbol{\mu}(k+1)=\mathbf{W}\boldsymbol{\mu}(k)-\alpha \mathbf{y}(k), \label{compact_mu} \\
	&\mathbf{x}(k+1)=\arg \min _{\mathbf{z}\in \mathcal {X}}\lbrace  f(\mathbf{z})-\boldsymbol{\mu}^{\intercal}(k+1) \mathbf{A}\mathbf{z} \rbrace, \label{compact_x} \\
	&\mathbf{y}(k+1)=\mathbf{W}\mathbf{y}(k)+\mathbf{A}\mathbf{x}(k+1)-\mathbf{A}\mathbf{x}(k). \label{deviation_tracking_1}
\end{align}

The following Lemma shows the fixed point of the proposed algorithm is consistent with the optimal solution of the problem (\ref{DRAP}) when Assumption 2 holds.

\emph{Lemma 1:} Under Assumption 2, any fixed point of the proposed algorithm satisfies the optimal conditions of constraint-coupled DRAP.

\emph{Proof:} It is not difficult to check that any fixed point $(\boldsymbol{\mu}_F, \mathbf{x}_F, \mathbf{y}_F)$ of the proposed algorithm satisfies
\begin{align}
	&(I_{nm}-\mathbf{W})\boldsymbol{\mu}_F=-\alpha \mathbf{y}_F, \label{fix_1}\\
	&\mathbf{x}_F=\arg \min _{\mathbf{z}\in \mathcal {X}}\lbrace  f(\mathbf{z})-\boldsymbol{\mu}_F^{\intercal}\mathbf{A}\mathbf{z} \rbrace, \label{fix_2} \\
	&(I_{nm}-\mathbf{W})\mathbf{y}_F=0. \label{fix_3}
\end{align}
Under Assumption 2, from (\ref{fix_1}) and (\ref{fix_3}), we have
\begin{align}
	(\mathbf{1}_n \otimes I_m)^\intercal \mathbf{y}_F =0, \ \ \mathbf{y}_F=\mathbf{1}_n \otimes y_F.
\end{align}
Thus, $\mathbf{y}_F=0$ holds and combining it with (\ref{fix_1}) yields 
\begin{align}
	&\boldsymbol{\mu }_F=\mathbf{1}_n \otimes \mu_F. \label{fix_4}
\end{align}
Moreover, with proper initialization, from (\ref{tracking_expression}) we have
\begin{align}
	(\mathbf{1}_n \otimes I_m)^\intercal (\mathbf{A}\mathbf{x}_F-\mathbf{d})=(\mathbf{1}_n \otimes I_m)^\intercal \mathbf{y}_F=0. \label{fix_5}
\end{align}
where (\ref{fix_2}), (\ref{fix_4}) and (\ref{fix_5}) correspond to  the optimal conditions in Proposition 1 which completes the proof. $\hfill\blacksquare$

With privacy concerns, agents exchange noise-masked information against potential attackers to prevent information disclosure. In \cite{Laplace_noise}, Laplace noise is shown to satisfy the necessary and sufficient condition of differential privacy. Furthermore, due to the noise accumulation in the tracking process as $\sum _{i=1}^n y_i=\sum _{i=1}^n A_ix_i - \sum _{i=1}^n d_i+\sum _{i=1}^n\sum_{t=0}^{n-1} \zeta_{i}(t)$, the added noise $\zeta_i$ needs to decay to ensure the convergence of $y_i$. Likewise, $\eta_i$ also needs to decay for the convergence of $\mu_i$ and $y_i$. Thus, the added noises are set to $\eta _{i}(k)\sim Lap(\theta _{ik}^{\eta })$, $\zeta _{i}(k)\sim Lap(\theta _{ik}^{\zeta })$ where $\theta _{ik}^{\eta }= d_{\eta_{i}}{q_i}^k, \theta _{ik}^{\zeta }= d_{\zeta_{i}}{q_i}^k, q_i\in(0,1)$ and let $\bar{q}=\max\{q_i\}$. That is, two diminishing Laplace noises are added to the communication process of $\mu_i$, $y_i$, respectively. Both of them are shown to be necessary in Section \ref{Differential Privacy}. 

Based on the mismatch tracking scheme and the information-masked protocol, a differentially private distributed mismatch tracking algorithm (diff-DMAC) is proposed which is shown in Algorithm 1 in details.

\begin{algorithm}[h]
	\caption{\textbf{: diff-DMAC Algorithm}}
	\label{algorithm}
	\begin{algorithmic}[1]
		\STATE {Initialization: $x_i(0)$ and $\mu_i(0)$ are arbitrarily assigned and $y_i(0)=A_ix_i(0)-d_i$.}
		\FOR{$k=0,1,...$}
		\STATE Adds noise $\eta_i(k), \zeta_i(k)$ to get $z_{\mu_i}(k)$ and $z_{y_i}(k)$ as \\
		$z_{\mu_j}(k) = \mu _j(k)+\eta _j(k)$, $z_{y_j}(k) = y_{j}(k)+\zeta_{j}(k)$.
		\STATE Broadcasts $z_{\mu_i}(k)$ and $z_{y_i}(k)$ to $j \in \mathcal{N}_i$.
		\STATE Receives $z_{\mu_j}(k)$ and $z_{y_j}(k)$ from $j \in \mathcal{N}_i$
		\STATE Updates $\mu_i(k+1)$ through \\
		$\mu _i(k+1)=\sum _{j=1}^{n}w_{ij}z_{\mu_j}(k)-\alpha y_i(k)$
		\STATE Updates $x_i(k+1)$ through \\
		$x_i(k+1)=\arg \min _{z\in \mathcal {X}_i}\lbrace f_i(z)-\mu _i^{\intercal}(k+1)A_iz \rbrace$
		\STATE Updates $y_i(k+1)$ through \\
		$y _i(k+1)=\sum _{j=1}^{n}w_{ij}z_{y_j}(k)+A_ix_i(k+1)-A_ix_i(k).$
		\ENDFOR
	\end{algorithmic}
\end{algorithm}
\section{Convergence Analysis for diff-DMAC}
\label{convergence analysis}
In this section, we will establish the linear convergence rate of the proposed algorithm in mean square. And the upper and lower bounds of the convergence accuracy is also provided. We set $d_i=0$, $\forall i$ in the following analysis for simplicity while nonzero scenarios can be analyzed similarly. 

Let $\nabla f(x)= [{\nabla f_1(x_1)}^\intercal, \cdots,{\nabla f_n(x_n)}^\intercal]^\intercal$ and $\mathbf{\bar{W}}=\frac{\mathbf{1}_n^{\intercal}\mathbf{1}_n}{n}\otimes I_m$, $\mathbf{\check{W}}=I_{nm}-\mathbf{\bar{W}}$, $\mathbf{\tilde{W}}= \mathbf{W} \mathbf{\check{W}}$. From relations (\ref{compact_mu}) and (\ref{deviation_tracking_1}), we have that 
\begin{align}
	&\bar{\boldsymbol{\mu }}(k+1)=\bar{\boldsymbol{\mu }}(k)+\bar{\boldsymbol{\eta }}(k)-\alpha\bar{\boldsymbol{y }}(k),  \label{mu_bar}\\
	&\bar{\boldsymbol{y }}(k+1)=\bar{\boldsymbol{y }}(k)+\bar{\boldsymbol{\zeta }}(k)+\mathbf{\bar{W}}\mathbf{A}({\mathbf{x}}(k+1)-{\mathbf{x}}(k)), \label{y_bar} \\
	&\check{\boldsymbol{\mu }}(k+1)=\mathbf{\tilde{W}}(\check{\boldsymbol{\mu }}(k)+\check{\boldsymbol{\eta }}(k))-\alpha\check{\boldsymbol{y }}(k), \label{mu_check} \\
	&\check{\boldsymbol{y }}(k+1)=\!\mathbf{\tilde{W}}(\check{\boldsymbol{y }}(k)\!+\!\!\check{\boldsymbol{\zeta }}(k))\!+\!\mathbf{\check{W}}\mathbf{A}(\mathbf{x(k\!+\!1)\!-\!x(k)}), \label{y_check} 
\end{align}
where $\bar{\boldsymbol{\mu }} = \mathbf{\bar{W}}{\boldsymbol{\mu }}$, $\check{\boldsymbol{\mu }} = \mathbf{\check{W}}{\boldsymbol{\mu }}$ and the doubly stochasticity of the weight matrix is used here.

\emph{Lemma 2:} For positive $a_1$, $a_2$, $a_3>0$, if $a_1a_2a_3>1$ holds, there exist $\gamma_1$, $\gamma_2$, $\gamma_3>0$ such that
\begin{align}
	\gamma_1< a_1 \gamma_2,\ \gamma_2< a_2 \gamma_3, \ \gamma_3< a_3 \gamma_1. \label{small_gain}
\end{align}
The proof is omitted. Note that for any $\beta>0$, $\gamma_i^\prime =\beta\gamma_i$, $i=1,2,3$, also satisfy the condition (\ref{small_gain}). Define $C = \{1+(\frac{\Vert \mathbf{A}\Vert^2\alpha^2}{\underline{\varphi}^2}-\frac{2\alpha}{\overline{L}})\cdot\underline{\lambda}(A A^\intercal)\}^{\frac{1}{2}}$. Based on Lemma 2, Theorem 1 is proved by induction. 

\emph{Theorem 1:} Given a constant $r \in (r_{LB},\infty)$, where $r_{LB}$ is defined as $r_{LB} = \max \{\bar{q}, C, \bar{\lambda}\}$. Suppose Assumptions 1,2 hold, $\alpha$ satisfies 
\begin{align}
	\alpha < \frac{\underline{\varphi}^2}{2\Vert \mathbf{A}\Vert^2 \overline{L}} \text{ and } \frac{(r-C)\underline{\varphi}}{\alpha \Vert \mathbf{A} \Vert}(\frac{(r-\bar{\lambda})^2\underline{\varphi}}{2\alpha \Vert \mathbf{A} \Vert}-1)>1. \label{alpha_1}
\end{align}
then the sequences generated by diff-DMAC satisfy
\begin{align}
	&\mathbb{E}[\Vert \check{\boldsymbol{\mu }}(k)\Vert] \leq \gamma_1r^k,
	\label{convergence_begin} \\
	&\mathbb{E}[\Vert \check{\boldsymbol{y}}(k)\Vert] \leq \gamma_2r^k, 
	\\
	&\mathbb{E}[\Vert \boldsymbol{x}(k)-\boldsymbol{x}^{\infty }\Vert] \leq \gamma_3r^k,  \label{x_k_x_infinity} 
	\\
	&\mathbb{E}[\Vert \bar{\boldsymbol{\mu }}(k)-\boldsymbol{\mu }^{\infty }\Vert] \leq \gamma_4r^k, \label{convergence_end} 
\end{align}
for some constant $\gamma_1, \gamma_2, \gamma_3, \gamma_4>0$, where $\mathbf{x}^\infty$ and $\boldsymbol{\mu}^\infty$ satisfy
\begin{align}
	&\bar{\mathbf{W }}\mathbf{A}{\mathbf{x}}^\infty=-\sum \limits_{k=1}^\infty\bar{\boldsymbol{\zeta }}(k), \ \bar{\mathbf{W }}{\boldsymbol{\mu }}^\infty={\boldsymbol{\mu }}^\infty, \label{infinity_1} \\
	&x_i^\infty=\arg \min _{z\in \mathcal {X}_i}\lbrace  f_i(z)-{\mu _i^\infty}^{\intercal} A_iz \rbrace , \forall i. \label{infinity_2} 
\end{align}

\emph{Proof:} Combining (\ref{mu_bar}) and (\ref{infinity_1}), we have that
\begin{align}
	\bar{\boldsymbol{\mu }}(k+1)-\bar{\boldsymbol{\mu }}^\infty=&\bar{\boldsymbol{\mu }}(k)-\bar{\boldsymbol{\mu }}^\infty-\alpha\bar{\mathbf{W }}\mathbf{A}({\mathbf{x}}(k)-{\mathbf{x}}^\infty) \notag \\
	&+\alpha \sum \limits_{t=k}^\infty \bar{\boldsymbol{\zeta }}(t)+\bar{\boldsymbol{\eta }}(k). \label{m4_1}
\end{align}

Under Assumption 1, the function 
\begin{align*}
	f_i^*(A_i^\intercal \mu _i) = \arg \max _{z\in \mathcal {X}_i}\lbrace \mu _i^{\intercal}A_iz-f_i(z) \rbrace
\end{align*} 
is differentiable and $f_i^*$ satisfies $\frac{1}{L_i}-$strongly convex and $\frac{1}{\varphi_i}-$Lipschitz smooth  \cite{conjugate_function}. Furthermore, it follows (\ref{compact_x}) and the Proposition B.25 in \cite{Danskins_theorem} that 
\begin{align}
	x_i(k)=\nabla f_i^*(A_i^\intercal \mu _i(k)), \forall i. \label{grad_f}
\end{align}

From (\ref{infinity_2}) and (\ref{grad_f}), we can obtain that
\begin{align}
	&\Vert \bar{\boldsymbol{\mu }}(k)-\boldsymbol{\mu }^{\infty } -\alpha \bar{\mathbf{W }} \mathbf{A} ({\mathbf{x}}(k)-{\mathbf{x}}^\infty)  \Vert \notag \\
	=&\Vert \bar{\boldsymbol{\mu }}(k)-\boldsymbol{\mu }^{\infty }-\alpha\bar{\mathbf{W }}\mathbf{A}(\nabla f^*(\mathbf{A} ^\intercal \boldsymbol{\mu}(k))-\nabla f^*(\mathbf{A}^\intercal \boldsymbol{\mu}^\infty)) \Vert \notag \\
	\leq& \left \{\Vert \bar{\boldsymbol{\mu }}(k)-\boldsymbol{\mu }^{\infty } \Vert^2+\alpha^2 \Vert \nabla f^*(\mathbf{A}^\intercal \bar{ \boldsymbol{\mu }}(k))-\nabla f^*(\mathbf{A}^\intercal \boldsymbol{\mu}^\infty) \Vert \right. \notag \\
	&  \left. -2\alpha(\bar{\boldsymbol{\mu }}(k)-\boldsymbol{\mu}^\infty)^\intercal\mathbf{A}(\nabla f^*(\mathbf{A}^\intercal\bar{\boldsymbol{\mu }}(k))-\nabla f^*(\mathbf{A}^\intercal \boldsymbol{\mu}^\infty)) \right \}^{\frac{1}{2}} \notag \\
	&  + \alpha \Vert \nabla f^*({\mathbf{A}^\intercal\boldsymbol{\mu }}(k))-\nabla f^*(\mathbf{A}^\intercal\bar {\boldsymbol{\mu}}(k)) \Vert \notag \\
	\leq& C \Vert \bar{\boldsymbol{\mu }}(k)-\boldsymbol{\mu }^{\infty } \Vert+\frac{\Vert \mathbf{A}\Vert^2 \alpha}{\underline{\varphi}} \Vert \check{\boldsymbol{\mu }}(k) \Vert, \label{m4_2}
\end{align}
where $C<1$ under (\ref{alpha_1}) and the last inequality follows from the strongly convexity and Lipschitz smoothness of $f^*$.

We will prove (\ref{convergence_begin})-(\ref{convergence_end}) by induction. When $\kappa = 0$, it is not difficult to find $\gamma_1, \gamma_2, \gamma_3, \gamma_4>0$ such that (\ref{convergence_begin})-(\ref{convergence_end}) hold. We assume that for all $\kappa \leq k$, (\ref{convergence_begin})-(\ref{convergence_end}) hold. Considering $\kappa=k+1$, we obtain from (\ref{m4_1}) and (\ref{m4_2}) that
\begin{align}
	&\mathbb{E}[\Vert \bar{\boldsymbol{\mu }}(k+1)-\boldsymbol{\mu }^{\infty } \Vert]  \notag \\
	&\leq C \mathbb{E}[\Vert \bar{\boldsymbol{\mu }}(k)-\boldsymbol{\mu }^{\infty } \Vert]+\frac{\Vert \mathbf{A}\Vert\alpha}{\underline{\varphi}} \mathbb{E}[\Vert \check{\boldsymbol{\mu }}(k) \Vert] \notag \\
	& \ \ \ +\mathbb{E}[ \Vert \alpha \sum \limits_{t=k}^\infty \bar{\boldsymbol{\zeta }}(t) \Vert] - \mathbb{E} [\Vert \bar{\boldsymbol{\eta }}(k) \Vert] \notag \\
	&\leq C\gamma_4r^k+\frac{\Vert \mathbf{A}\Vert \alpha}{\underline{\varphi}}\gamma_1r^k+nm\bar{q}^k\bar{d}_{\eta }+\frac{\alpha nm\bar{q}^k\bar{d}_{\zeta }}{1-\bar{q}} . \label{m4_final}
\end{align}

Combining (\ref{mu_check})-(\ref{y_check}) gives
\begin{align}
	\mathbb{E}[\Vert \check{\boldsymbol{\mu }}(k+1) \Vert] 
	&\leq \bar{\lambda}\mathbb{E}[\Vert \check{\boldsymbol{\mu }}(k)\Vert] +\alpha\mathbb{E}[\Vert \check{\boldsymbol{y}}(k) \Vert] +\bar{\lambda}\mathbb{E}[\Vert \check{\boldsymbol{\eta }}(k)\Vert] \notag \\
	&\leq \bar{\lambda}\gamma_1r^k+\alpha \gamma_2 r^k +\bar{\lambda}nm\bar{q}^k\bar{d}_{\eta }, \label{m1_final}
\end{align}
\begin{align}
	\mathbb{E}[\Vert \check{\boldsymbol{y}}(k+1) \Vert] 
	&\leq \bar{\lambda}\mathbb{E}[\Vert \check{\boldsymbol{y}}(k)\Vert] +\Vert \mathbf{A}\Vert \mathbb{E}[\Vert \boldsymbol{x}(k) -\boldsymbol{x}^\infty \Vert] \notag \\
	& \ \ \ +\Vert \mathbf{A}\Vert\mathbb{E}[\Vert \boldsymbol{x}(k+1) -\boldsymbol{x}^\infty \Vert] +\bar{\lambda}\mathbb{E}[\Vert \check{\boldsymbol{\zeta }}(k)\Vert] \notag \\
	&\leq \bar{\lambda}\gamma_2r^k+\Vert \mathbf{A}\Vert \gamma_3 r^k +\bar{\lambda}nm\bar{q}^k\bar{d}_{\zeta } \notag \\
	& \ \ \ +\Vert \mathbf{A}\Vert\mathbb{E}[\Vert \boldsymbol{x}(k+1) -\boldsymbol{x}^\infty \Vert]. \label{m2_final}
\end{align}

Under Assumption 1, (\ref{grad_f}) implies that
\begin{align}
	\mathbb{E}[\Vert \boldsymbol{x}(k+1) -\boldsymbol{x}^\infty \Vert]=\mathbb{E}[\Vert \nabla f^*(\mathbf{A}\boldsymbol{\mu}(k+1))- \nabla f^*(\mathbf{A}\boldsymbol{\mu}^\infty) \Vert] \notag \\
	\leq\frac{\Vert \mathbf{A}\Vert}{\underline{\varphi}} \mathbb{E}[\Vert  \check{\boldsymbol{\mu}}(k+1) \Vert] + \frac{\Vert \mathbf{A}\Vert}{\underline{\varphi}} \mathbb{E}[\Vert \bar{\boldsymbol{\mu}}(k+1)-\boldsymbol{\mu}^\infty \Vert], \label{m3_final}
\end{align}
where the last inequality is due to the Lipschitz smoothness of $f^*$.

Note that constant $\bar{\lambda}, C \in(0,1)$, due to Lemma 2, there exist $\gamma_1, \gamma_2, \gamma_3, \gamma_4 >0 $ such that 
\begin{align}
	C\gamma_4+\frac{\Vert \mathbf{A}\Vert \alpha}{\underline{\varphi}}\gamma_1<\gamma_4r, \label{m_begin} \\
	\bar{\lambda}\gamma_1+\alpha \gamma_2<\gamma_1r,  \\
	\bar{\lambda}\gamma_2+2\Vert \mathbf{A}\Vert \gamma_3<\gamma_2r,  \\
	\frac{\Vert \mathbf{A}\Vert}{\underline{\varphi}}(\gamma_1+\gamma_4)<\gamma_3.\label{m_end}
\end{align}
Basd on equations (\ref{m_begin})-(\ref{m_end}), taking suffciently large $\gamma_1, \gamma_2, \gamma_3, \gamma_4 >0 $ yields
\begin{align}
	C\gamma_4r^k+\frac{\Vert \mathbf{A}\Vert \alpha}{\bar{\varphi}}\gamma_1r^k+nm\bar{q}^k\bar{d}_{\eta }+\frac{\alpha nm\bar{q}^k\bar{d}_{\zeta }}{1-\bar{q}}& \leq \gamma_4r^{k+1}, \notag \\
	\bar{\lambda}\gamma_1r^k+\alpha \gamma_2 r^k +\bar{\lambda}nm\bar{q}^k\bar{d}_{\eta }& \leq \gamma_1r^{k+1}, \notag \\
	\bar{\lambda}\gamma_2r^k+2\Vert \mathbf{A}\Vert \gamma_3r^k& \leq \gamma_2r^{k+1}, \notag \\
	\frac{\Vert \mathbf{A}\Vert}{\underline{\varphi}}(\gamma_1+\gamma_4)r^{k+1}& \leq \gamma_3r^{k+1}, \notag
\end{align}
where we use the definition of $r_{LB}$ and substituting these equations into (\ref{m4_final})-(\ref{m3_final}) completes the proof. $\hfill\blacksquare$

Based on Theorem 1, we further prove the linear convergence property of diff-DMAC and quantify its convergence accuracy in Theorem 2.

\emph{Theorem 2:} Suppose the conditions in Theorem 1 are satisfied. If $\alpha$ further satisfies
\begin{align}
	\alpha<\frac{\underline{\varphi}[-(1\!-\!C)\!+\!\sqrt{(1\!-\!C)^2\!+\!2(1\!-\!C)(1\!-\!\bar{\lambda})^2}]}{2\Vert \mathbf{A}\Vert}, \label{alpha_4}
\end{align}
the sequence $\{\mathbf{x(k)}\}_{k\geq 0}$ generated by diff-DMAC will converge linearly to the neighbor of the optimum $\mathbf{x}^\star$ in mean square and we have 
\begin{align}
	\frac{1}{n^2\Vert \mathbf{A}\Vert^2}N_\zeta \leq \mathbb {E}[ \Vert \boldsymbol{x}^{\infty }-\boldsymbol{x}^\star\Vert ^2] \leq \frac{\overline{L}^2}{n\underline{\varphi }^2 \underline{\lambda}(\mathbf{A} \mathbf{A}^\intercal)}N_\zeta , \label{expectation_error}
\end{align}
where $N_\zeta $ is defined as $N_\zeta = \sum _{i=1}^n\frac{2md_{\zeta i}^2}{1-q_i^2}$.

\emph{Proof:} Note that $\bar{\lambda}, C \in(0,1)$, it is not difficult to see that with $\alpha$ chosen to satisfy (\ref{alpha_4}), we can always find a constant $r \in (r_{LB},1)$ such that equations (\ref{convergence_begin})-(\ref{convergence_end}) hold. Thus, the convergence of diff-DMAC is proved by Theorem 1.

Due to (\ref{infinity_1}) and $\bar{\mathbf{W}}\mathbf{A}{\mathbf{x}}^\star=0$, we have
\begin{align}
	\mathbb{E}[\Vert \bar{\mathbf{W}}\mathbf{A}({\mathbf{x}}^\infty - {\mathbf{x}}^\star) \Vert^2]=\mathbb{E}[\Vert \sum\limits_{t=0}^\infty \bar{\boldsymbol{\zeta }}(t)\Vert^2]=\frac{1}{n^2}\sum \limits_{i=1}^n \frac{2md_{\zeta_i}^2}{1-q_{\zeta_i}^2}. \label{opt_value}
\end{align}

Since $\mathbf{\bar{W}}=\frac{\mathbf{1}^{\intercal}\mathbf{1}}{n}\otimes I_m$, it is easy to obtain that
\begin{align}
	\mathbb{E}[\Vert \bar{\mathbf{W}}\mathbf{A}({\mathbf{x}}^\infty - {\mathbf{x}}^\star) \Vert^2] \leq \Vert \mathbf{A}\Vert^2\mathbb{E}[\Vert \mathbf{x}^\infty -\mathbf{x}^\star \Vert ^2]. \label{lower_bound}
\end{align}

The optimality condition of the $\mathbf{x}-$update (\ref{compact_x}) implies
\begin{align}
	(\nabla f_i(x_i^\infty)-A_i ^\intercal \mu_i^\infty)^\intercal(x_i^\star-x_i^\infty) &\geq 0 \label{xupdate_opt_1}, \\
	(\nabla f_i(x_i^\star)-A_i^\intercal\mu_i^\star)^\intercal(x_i^\infty-x_i^\star) &\geq 0 \label{xupdate_opt_2}, \forall i.
\end{align}

Based on (\ref{xupdate_opt_1}) and (\ref{xupdate_opt_2}), we obtain
\begin{align}
	&(\mu_i^\infty-\mu_i^\star)^\intercal A_i(x_i^\infty-x_i^\star) \notag \\
	\geq& (\nabla f_i( x_i^\infty)-\nabla f_i( x_i^\star))^\intercal(x_i^\infty-x_i^\star) \notag \\
	\geq& \frac{\underline{\varphi}}{\bar{L}} \Vert A_i^\intercal \mu_i^\infty-A_i^\intercal \mu_i^\star \Vert \Vert x_i^\infty-x_i^\star\Vert, \forall i \in \mathcal{N}. \label{xupdate_opt_3}
\end{align}

The facts that $\bar{\mathbf{W }}{\boldsymbol{\mu }}^\infty={\boldsymbol{\mu }}^\infty$ and $\bar{\mathbf{W }}{\boldsymbol{\mu }}^\star={\boldsymbol{\mu }}^\star$ imply 
\begin{align}
	&\mu_i^\infty-\mu_i^\star= \frac{1}{n}\sum \limits_{i=1}^n \mu_i^\infty-\frac{1}{n}\sum \limits_{i=1}^n \mu_i^\star, \ \forall i \in \mathcal{V} \label{xupdate_opt_4}, \\
	&(\boldsymbol{\mu}^\infty-\boldsymbol{\mu}^\star)^\intercal \mathbf{A} (\mathbf{x}^\infty-\mathbf{x}^\star) \notag\\ =&(\boldsymbol{\mu}^\infty-\boldsymbol{\mu}^\star)^\intercal \bar{\mathbf{W }} \mathbf{A} ({\mathbf{x}}^\infty-{\mathbf{x}}^\star)\notag \\
	\leq& \Vert \boldsymbol{\mu}^\infty-\boldsymbol{\mu}^\star\Vert \Vert \bar{\mathbf{W }} \mathbf{A} ({\mathbf{x}}^\infty-{\mathbf{x}}^\star)\Vert, \label{xupdate_opt_5}
\end{align}

Then, combining (\ref{xupdate_opt_3})-(\ref{xupdate_opt_5}) yields
\begin{align}
	\mathbb{E}[\Vert \mathbf{x}^\infty -\mathbf{x}^\star\Vert ^2] \leq \frac{\bar{L}^2}{n\underline{\varphi}^2 \underline{\lambda}(\mathbf{A} \mathbf{A}^\intercal)}\mathbb{E}[\Vert \bar{\mathbf{W }} \mathbf{A} ({\mathbf{x}}^\infty-{\mathbf{x}}^\star)\Vert ^2]. \label{upper_bound}
\end{align}

Substituting (\ref{opt_value}) into (\ref{lower_bound}), (\ref{upper_bound}) completes the  proof. $\hfill\blacksquare$

Therefore, the convergence error is independent of initial conditions of the algorithm and it only related to noise $\boldsymbol{\zeta }$. 

It follows from relations (\ref{x_k_x_infinity}) and (\ref{expectation_error}) that 
\begin{align}
	&\mathbb {E}\lbrace \Vert \boldsymbol{x}(k)-\boldsymbol{x}^\star\Vert ^2\rbrace  \notag \\
	\leq & 2\lbrace\mathbb {E} [\Vert \boldsymbol{x}(k)-\boldsymbol{x}^\infty\Vert ^2]+\mathbb {E} [\Vert \boldsymbol{x}^{\infty }-\boldsymbol{x}^\star\Vert ^2]\rbrace \notag \\
	\leq & 2[\gamma_3r^k+\frac{\overline{L}^2}{n\underline{\varphi }^2 \underline{\lambda}(\mathbf{A} \mathbf{A}^\intercal)}N_\zeta]. \label{complexity}
\end{align}
Define $C_\zeta=\frac{\overline{L}^2}{n\underline{\varphi }^2 \underline{\lambda}(\mathbf{A} \mathbf{A}^\intercal)}N_\zeta$. For accuracy $\varepsilon>C_\zeta$, it can be derived from relation (\ref{complexity}) that the complexity of the proposed algorithm is $O(\ln(\frac{1}{\varepsilon-C_\zeta}))$.

\section{Differential Privacy}
\label{Differential Privacy}
In this section, we present the privacy preserving performance  of the proposed algorithm by theoretical analysis.

\emph{Theorem 3:} Suppose the conditions in Theorem 2 hold, if $q_{i_0}$ satisfies \begin{align}
	q_{i_0}\in (\frac{\alpha \Vert A_{i_0} \Vert ^2 +\Vert A_{i_0} \Vert \sqrt{\alpha ^2\Vert A_{i_0} \Vert ^2+4\alpha \varphi _{i_0}}}{2\varphi _{i_0}},1), \label{q_bound}
\end{align} 
the proposed diff-DMAC preserves $\epsilon_{i_0}$-differential privacy for the cost function of agent $i_0$ , where
\begin{align*}
	\epsilon_{i_0}=(\frac{1}{\alpha d_{\zeta_{i_0}}}+\frac{1}{d_{\eta_{i_0}}})\frac{\alpha \varphi_{i_0} \delta \Vert A_{i_0} \Vert}{\varphi _{i_0}q_{i_0}^2-\alpha \Vert A_{i_0} \Vert ^2 q_{i_0}-\alpha \Vert A_{i_0} \Vert ^2}. 
\end{align*} 

\emph{Proof:} Two function sets have the same initialization and all auxiliary variables are known to attackers as mentioned in Section \ref{Problem Formulation}. Therefore, exchanged variables must be same, i.e., $\boldsymbol{z}_{\boldsymbol{\mu }}^{(1)}(k)=\boldsymbol{z}_{\boldsymbol{\mu }}^{(2)}(k)$ and $\boldsymbol{z}_{\boldsymbol{y}}^{(1)}(k)=\boldsymbol{z}_{\boldsymbol{y}}^{(2)}(k)$, otherwise, the two function sets will be easily distinguished. Therefore, for any $i\ne i_0$, due to the realation $f_i^{(1)}=f_i^{(2)}$, the added noise should satisfy
\begin{equation*}
	\zeta _i^{(1)}(k)=\zeta _i^{(2)}(k), \ \eta _i^{(1)}(k)=\eta _i^{(2)}(k), \forall k, \forall i\ne i_0. 
\end{equation*}
And for $i = i_0$, since $f_i^{(1)} \ne f_i^{(2)}$ the noise should satisfy
\begin{align}
	\Delta \eta _{i_0}(k)=-\Delta \mu _{i_0}(k), \ \Delta \zeta _{i_0}(k)=-\Delta y_{i_0}(k), \label{delta_1}
\end{align}
where $\Delta \zeta _{i_0}(k)=\zeta _{i_0}^{(1)}(k)-\zeta _{i_0}^{(2)}(k)$, $\Delta \eta _{i_0}(k)=\eta _{i_0}^{(1)}(k)-\eta _{i_0}^{(2)}(k)$, $\Delta \mu _{i_0}(k)=\mu _{i_0}^{(1)}(k)-\mu _{i_0}^{(2)}(k)$ and $\Delta y_{i_0}(k)=y_{i_0}^{(1)}(k)-y_{i_0}^{(2)}(k)$. Furthermore, it follows that
\begin{align}
	\Delta \mu _{i_0}(k+1)=-\alpha \Delta y_{i_0}(k), \ \Delta y_{i_0}(k+1)= A_{i_0}\Delta g_{i_0}(k),  \label{delta_2}
\end{align}
where $\Delta g_{i_0}(k)=\Delta x_{i_0}(k+1)-\Delta x_{i_0}(k)$ and $\Delta x_{i_0}(k)=x_{i_0}^{(1)}(k)-x_{i_0}^{(2)}(k)$.

Following the analysis in \cite{TAC_Ding}, we need to obtain the upper bound of $\frac{f_{\zeta \eta }(\boldsymbol{\zeta }^{(1)}, \boldsymbol{\eta }^{(1)})}{f_{\zeta \eta }(\mathbb{B}(\boldsymbol{\zeta }^{(1)}, \boldsymbol{\eta }^{(1)}))}$ to quantify the privacy level $\epsilon_{i_0}$,
\begin{align}
	&\frac{\mathbb {P}\lbrace \boldsymbol{\zeta }, \boldsymbol{\eta }\in \Omega |Z_{T^{(1)}}(\boldsymbol{\zeta }, \boldsymbol{\eta })\in \mathcal {O}\rbrace }{\mathbb {P}\lbrace \boldsymbol{\zeta }, \boldsymbol{\eta }\in \Omega |Z_{T^{(2)}}(\boldsymbol{\zeta }, \boldsymbol{\eta })\in \mathcal {O}\rbrace } \label{sup_definition} \leq \sup_{(\boldsymbol{\zeta },\boldsymbol{\eta })} \frac{f_{\zeta \eta }(\boldsymbol{\zeta }^{(1)}, \boldsymbol{\eta }^{(1)})}{f_{\zeta \eta }(\mathbb{B}(\boldsymbol{\zeta }^{(1)}, \boldsymbol{\eta }^{(1)}))},
\end{align}
where $Z_{T^{(l)}}(\boldsymbol{\zeta }^{(l)}  \boldsymbol{\eta }^{(l)})=\{\boldsymbol{z}_{\boldsymbol{\mu}}^{(l)}(k), \boldsymbol{z}_{\boldsymbol{y}}^{(l)}(k)\}_{k\geq0}, l=1, 2$, belongs to attacker's observation.

For $l\in \{1,2\}$, the $\mathbf{x}$-update can be rewritten as
\begin{align*}
	x_{i_0}^{(l)}(k+1)=\arg \min _{z\in \mathcal {X}_{i_0}^{(l)}}\lbrace  f_{i_0}^{(l)}(z)-{\mu _{i_0}^{(l)}}^{\intercal}(k+1)A_{i_0}z \rbrace. 
\end{align*}

From Definition 1, we have $x_{i_0}^{(1)}, x_{i_0}^{(2)}-\delta' \in \mathcal {X}_{i_0}^{(1)}$ and  $x_{i_0}^{(2)}, x_{i_0}^{(1)}+\delta' \in \mathcal {X}_{i_0}^{(1)}$. Similar to the analysis of Theorem 2, the optimal condition of the written $\mathbf{x}$-update yields
\begin{align}
	&(\mu_{i_0}^{(1)}(k)-\mu_{i_0}^{(2)}(k))^\intercal A_{i_0}(x_{i_0}^{(1)}(k)+\delta'-x_{i_0}^{(2)}(k)) \notag \\
	&\geq \frac{\varphi_{i_0}}{\Vert A_{i_0} \Vert^2} \Vert  A_{i_0}(x_{i_0}^{(1)}(k)+\delta'-x_{i_0}^{(2)}(k)) \Vert^2. \label{diff_x_opt}
\end{align}
The last inequality uses the relation $ \nabla f_{i_0}^{(1)}(x_{i_0}^{(1)}(k))=\nabla f_{i_0}^{(2)}(x_{i_0}^{(1)}(k)+\delta ^{\prime })$. Then applying Cauchy–Schwarz inequality to the LHS of (\ref{diff_x_opt}) yields
\begin{align}
	& \Vert A_{i_0}(x_{i_0}^{(1)}(k)+\delta'-x_{i_0}^{(2)}(k)) \Vert  \leq \frac{\Vert A_{i_0} \Vert^2}{\varphi_{i_0}} \Vert \Delta\mu_{i_0}(k) \Vert \label{diff_x_to_mu},
\end{align}
 where the Lipschitz smoothness is used.
 
From (\ref{delta_1}), (\ref{delta_2}) and the relation (\ref{diff_x_to_mu}), we have 
\begin{align*}
	\left\Vert \Delta \eta _{i_0}(k+1)\right\Vert =&\alpha\left\Vert A_{i_0}\Delta g_{i_0}(k-1)\right\Vert _2 \notag \\
	=&\alpha \Big \Vert A_{i_0}(x_{i_0}^{(1)}(k)-x_{i_0}^{(2)}(k)+\delta^\prime) \notag \\
	&-A_{i_0}\left( x_{i_0}^{(1)}(k-1)-x_{i_0}^{(2)}(k-1)+\delta^\prime \right)\Big \Vert  \notag \\
	\leq &\frac{\alpha \Vert A_{i_0} \Vert^2}{\varphi _{i_0}}\left(\left\Vert \Delta \eta _{i_0}(k)\right\Vert _2+\left\Vert \Delta \eta _{i_0}(k-1)\right\Vert \right).
\end{align*}
The two roots of $\varphi _{i_0}x^2-\alpha \Vert A_{i_0} \Vert^2 x-\alpha\Vert A_{i_0} \Vert^2 =0$, denoted as $\tau_1, \tau_2$ are different and all lie in $(-1,1)$.
Therefore, from the above relation, we have
\begin{align}
	\Vert \Delta \eta _{i_0}(k)\Vert _2 &\leq \frac{\alpha \delta \Vert A_{i_0}\Vert}{\tau _1-\tau _2}(\tau _1^{k-1}-\tau _2^{k-1}), \label{dengbishulie}
\end{align}
since $\Vert \Delta \eta _{i_0}(1)\Vert =0$ and $\Vert \Delta \eta _{i_0}(2)\Vert < \alpha \delta \Vert A_{i_0} \Vert $ which can be derived from (\ref{delta_1}) and (\ref{delta_2}). These two equations also give that $\alpha |\Delta \zeta _{i_0}(k)|=|\Delta \eta _{i_0}(k+1)|, \forall k$. Thus, 
\begin{align}
	\frac{f_{\zeta \eta }(\boldsymbol{\zeta }^{(1)}, \boldsymbol{\eta }^{(1)})}{f_{\zeta \eta }\left(\mathbb{B}(\boldsymbol{\zeta }^{(1)}, \boldsymbol{\eta }^{(1)})\right)}&=\prod _{k=1}^{\infty }\frac{f_L(\zeta _{i_0}^{(1)}(k),\theta _{ik}^{\zeta })f_L(\eta _{i_0}^{(1)}(k),\theta _{ik}^{\eta })}{f_L(\zeta _{i_0}^{(2)}(k),\theta _{ik}^{\zeta })f_L(\eta _{i_0}^{(2)}(k),\theta _{ik}^{\eta })}  \notag \\
	&\leq \prod _{k=1}^{\infty }e^{\frac{|\Delta \zeta _{i_0}(k)|}{d_{\zeta i_0}q_{i_0}^k}}\prod _{k=2}^{\infty }e^{\frac{|\Delta \eta _{i_0}(k)|}{d_{\eta i_0}q_{i_0}^k}}  \notag\\
	&\leq e^{\left(\frac{1}{\alpha d_{\zeta i_0}}+\frac{1}{d_{\eta i_0}}\right)\frac{\alpha \varphi _{i_0}\delta \Vert A_{i_0} \Vert}{D }}\notag \\ 
	&= e^{\epsilon _{i_0}}, \label{diff_result}
\end{align}
where $D=\varphi _{i_0}q_{i_0}^2-\alpha \Vert A_{i_0} \Vert ^2 q_{i_0}-\alpha \Vert A_{i_0} \Vert ^2$ and the last inequality follows the the sum formula of geometric series. The inequality $D >0$ holds under (\ref{q_bound}).

Combining (\ref{sup_definition}) and (\ref{diff_result}) yields the result corresponding to the definition of differential privacy of the cost function in Definition 2. Thus, the proof completes. $\hfill\blacksquare$

\emph{Remark 1:}
The privacy level $\epsilon$ in Theorem 3 is related to both noise $\boldsymbol{\zeta}$ and $\boldsymbol{\eta}$. Both of them are necessary, since if at least one of them are set to zero, there does not exist a finite number $\epsilon$, i.e., the differential privacy cannot be preserved.

Since the convergence accuracy is only related to the noise $\boldsymbol{\zeta }$, the optimal value of $\epsilon_i$ can be obtained by setting $d_{\eta_i} \rightarrow \infty$ and it gives 
\begin{align*}
	\epsilon _{i}^\star=\frac{\varphi _{i}\delta \Vert A_{i}\Vert }{d_{\zeta i}(\varphi _{i}q_{i}^2-\alpha \Vert A_{i} \Vert ^2 q_{i}-\alpha \Vert A_{i} \Vert ^2)},\forall i \in \lbrace 1,\cdots ,n \rbrace.
\end{align*}
Combining the above one with (\ref{expectation_error}) shows that larger $d_{\zeta_{i}}$ and $q_i$ lead to a lower $\epsilon_{i}^\star$. However, the convergence error can be larger at the same time. That is, big and slow diminishing added noise leads to a high level privacy but a low accuracy which shows the trade-off between accuracy and privacy.

\section{Numerical Experiments} 
\label{Numerical Experiments}
In this section, we conduct numerical experiments to verify
our theoretical analysis. The proposed algorithm is tested on a multi-MG system with 14 microgrids \cite{multi_microgrid_systems}. Regarding the communication network, we generate an undirected connected graph by adding random links to a ring network. We consider the following problem with quadratic cost functions $f_i(x_i)=u_ix_i^2+v_ix_i+w_i, \forall i$
\begin{align}
	\min _{\boldsymbol{x}\in \mathbb {R}^{n \times 1}}f(\boldsymbol{x})&=\sum _{i=1}^{n}f_i(x_i) \notag \\
	s.t. \sum _{i=1}^{n} a_i x_i=\sum _{i=1}^{n}d_i , &\  x_i^{\min} \leq x_i \leq x_i^{\max},\forall i. \label{simulation}
\end{align} 
where coefficients $a_i$ are randomly chosen and the parameters
of the generators are adopted from \cite{economic_dispatch}.

\begin{figure}[thpb]
	\centering
	\subfigure[]{
		\includegraphics[width=4cm]{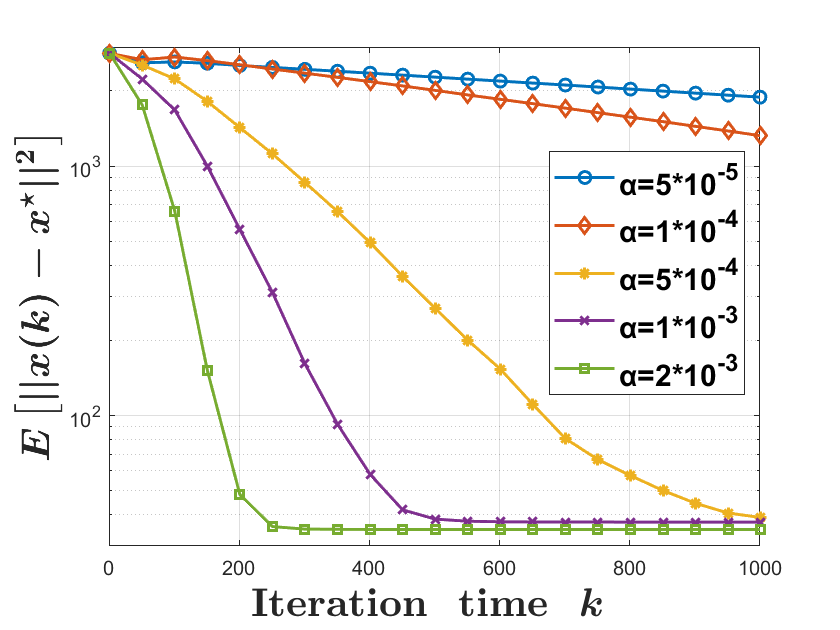}
	}
	\subfigure[]{
		\includegraphics[width=4cm]{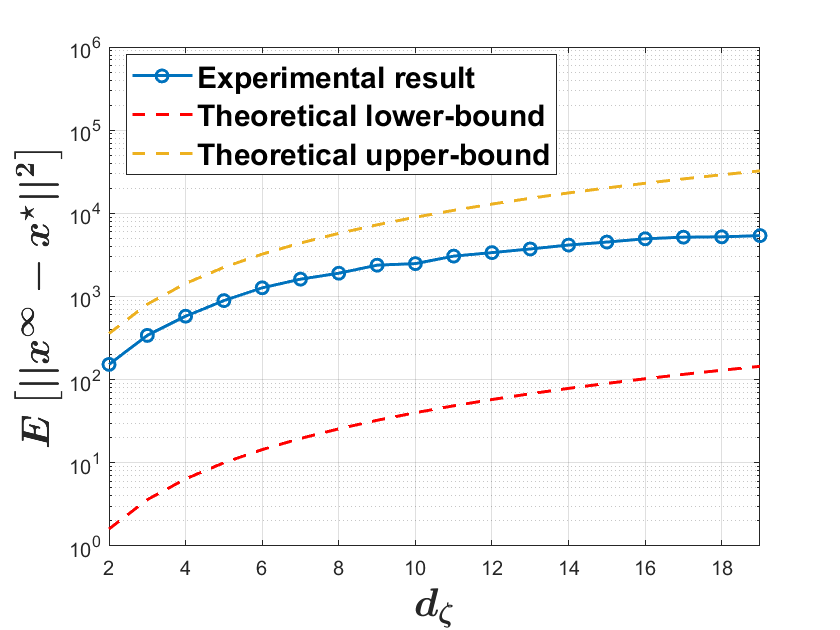}
	}
	\caption{(a) MSE versus Iteration time $k$ with different stepsizes; (b) Relation between MSE and $d_\zeta.$}
	\label{MSE_alpha_d}
\end{figure}

Set total load demand $d_{total}=231MW$, the exchanged information is masked by independent Laplace noise with $q=0.98, d_\eta=d_\zeta=1$ and we apply diff-DMAC to problem (\ref{simulation}). Fig. \ref{MSE_alpha_d} (a) shows how the mean square error (MSE) $\mathbb{E}[\Vert \mathbf{x}(k)-\mathbf{x}^\star \Vert^2]$ changes with iteration time $k$, where the expected errors are approximated by averaging over 100 simulation results. The results validate the linear convergence rate of the proposed algorithm. Fix $q=0.98$, the relation between error and $d_\zeta$ is shown in Fig. \ref{MSE_alpha_d} (b). We find that as $d_\zeta$ becomes larger, the MSE of $\mathbf{x}$ increases, which means
larger noise brings lower accuracy. Moreover, the experimental result is strictly between the lower bound and upper bound given in Theorem 2.

\begin{figure}[thpb]
	\centering
	\subfigure[]{
		\includegraphics[width=4cm]{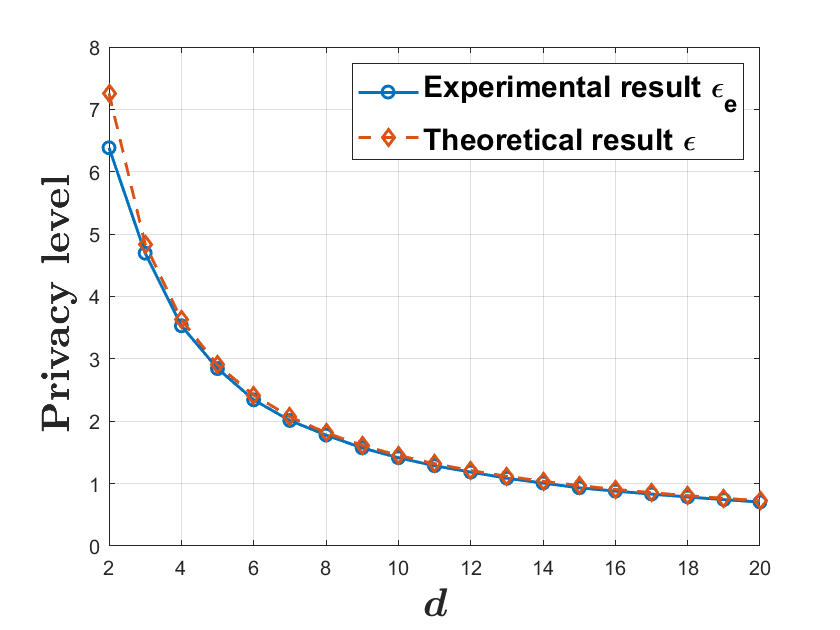}
	}
	\subfigure[]{
		\includegraphics[width=4cm]{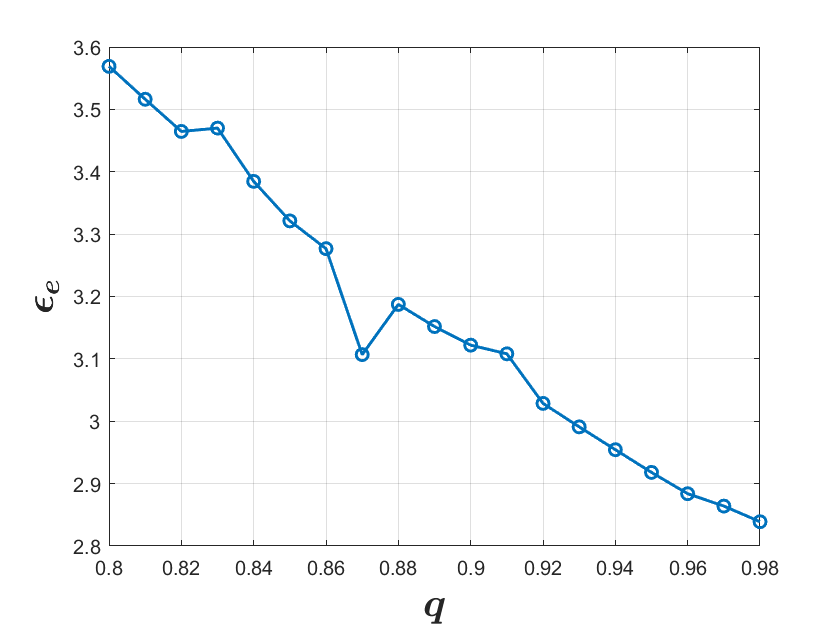}
	}
	\caption{(a) Privacy level versus $d_\zeta$; (b) $\epsilon_e$ versus decay coefficient $q$.}
	\label{privacy_validate}
\end{figure}

\begin{figure}[thpb]
	\centering
	\includegraphics[width=4cm]{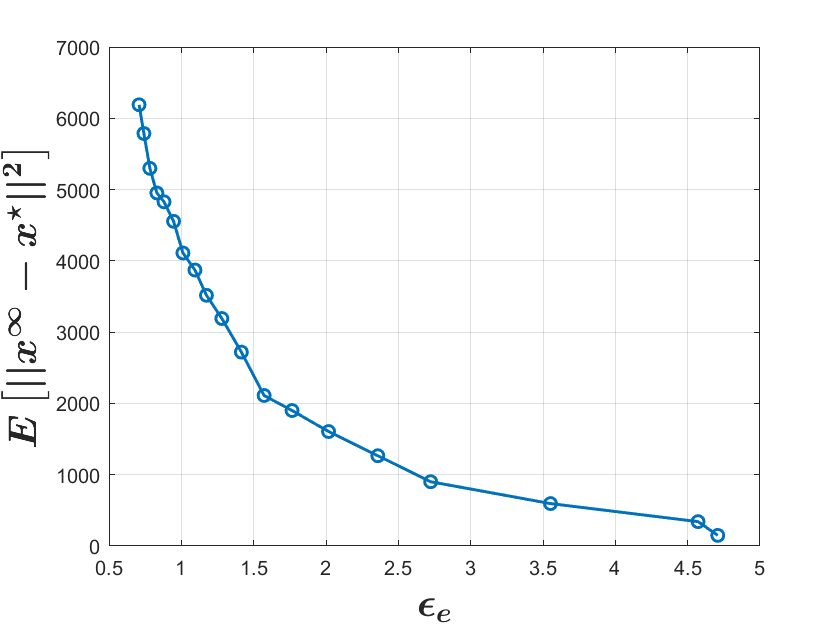}
	\caption{Trade-off between convergence accuracy and privacy level.}
	\label{trade-off}
\end{figure}

As for the differential privacy of diff-DMAC, Theorem 3 is numerically tested with different $d_\zeta$ as shown in Fig. \ref{privacy_validate}(a). The experimental result is upper bounded by  the theoretical result, i.e. $\epsilon_e\leq \epsilon$, which verifies Theorem 3. Combining Fig. \ref{privacy_validate}(a) with Fig. \ref{privacy_validate}(b) shows that a bigger and slower decaying noise leads to a better privacy level. However, large noise also influences the convergence accuracy and Fig. \ref{trade-off} validates the aforementioned trade-off between privacy and accuracy.

\section{Conclusion} 
In this paper, a differentially private distributed mismatch tracking algorithm has been proposed to solve constraint-coupled resource allocation problems with privacy concerns. Its linear convergence property has been established for strongly convex and Lipschitz-smooth cost functions. Then, the differential privacy of the proposed algorithm is theoretically proved and we also characterize the trade-off between accuracy and privacy. The theoretical results have been examined by numerical experiments.
\label{Conclusion}

\bibliographystyle{IEEEtran}
\bibliography{ref}

\end{document}